\documentclass{article}\usepackage{amsmath,amssymb,latexsym,amsfonts,enumerate,theorem,longtable,amscd}
\newtheorem{definition}{Definition}[section]\newtheorem{theorem}{Theorem}[section]
\newtheorem{corollary}{Corollary}[section]\newtheorem{lemma}{Lemma}[section]
\begin{document}

\begin{center}
{\bf\large Some topics on the $Z_3$-graded exterior algebras}
\end{center}

\begin{center}
Salih Celik\footnote{E-mail: sacelik@yildiz.edu.tr} and Sultan A. Celik

Department of Mathematics, Yildiz Technical University, DAVUTPASA-Esenler, Istanbul, 34210 TURKEY.
\end{center}



\begin{abstract}
A $Z_3$-graded Hopf algebra structure of exterior algebra with two generators is introduced. Two covariant differential calculus on the $Z_3$-graded exterior algebra are presented. Using the generators and their partial derivatives a Grassmann-Heisenberg algebra is constructed. An R-matrix which satisfies graded Yang-Baxter equations is obtained. A $Z_3$-graded universal enveloping algebra $U_q(\widetilde{gl}(2))$ is constructed with the quadratic elements of the Grassmann-Weyl algebra.
\end{abstract}

\section{Introduction}\label{sec1}

It is well known that a "deformation" of a mathematical structure is a family of similar mathematical structures, which depends on a parameter (or parameters) such that the original structure corresponds to a particular chosen initial value for the parameter.

Quantum plane \cite{Manin1} is a well known example in quantum group theory. The natural extension to $Z_2$-graded space was introduced in \cite{Manin2}. In the last decade, the $Z_2$-graded structures has been extensively developed to the $Z_3$-graded structures (see, for example, \cite{Ab1}-\cite{Celik2}, \cite{Chung}, \cite{Kerner1}), and $Z_N$-graded constructions \cite{Coq}, \cite{Dub}, \cite{Kapr}, \cite{KerNie}.

The quantum plane, on the contrary, is a more peculiar object having the unusual property of not being commutative. The noncommutativity of these variables provides an example of noncommutative geometry \cite{Con}. It can be obtained from \cite{Woro}, that one can set up a consistent differential calculus on the noncommutative spaces. Quantum (hyper)plane relations have been a model to establish a noncommutative differential calculus \cite{WZ}. This calculus is covariant under the action of quantum group. In that approach, differential forms are defined in terms of $q$-commuting coordinates, and the differential and algebraic properties of quantum groups acting on these spaces are obtained from the properties of the spaces. The present work starts a $Z_3$-graded version of the exterior plane, denoted by $\tilde{\mathbb R}_q^{0|2}$, where $q$ is a cubic root of unity. In this case, of course, it will not go back to the original objects. A $Z_3$-graded differential calculus on a space involves functions on the space, first and second order differential of coordinates. Another important property of such a calculus is that the operator {\sf d} satisfies ${\sf d}^3:={\sf d}\circ{\sf d}\circ{\sf d}=0$ instead of ${\sf d}^2=0$ such that ${\sf d}^2\ne0$. In Sect. 4, we set up differential calculi on the $Z_3$-graded exterior plane. Boht of them are covariant under the action of the $Z_3$-graded (quantum) group. In Sect. 6, using the generators and their partial derivatives we construct a Grassmann-Heisenberg algebra. In Sect. 7, we introduce an $R$-matrix preserving all noncommutative structures. We introduce a $Z_3$-graded universal enveloping algebra $U_q(\widetilde{gl}(2))$ with the quadratic elements of the Grassmann-Weyl algebra in Sect. 8.

\section{$Z_3$-graded exterior algebras}\label{sec2}

The aim of this section is to introduce the $Z_3$-graded version of the exterior plane and its dual. It is known that the Manin's quantum plane is introduced as a $q$-deformation of commutative plane in the sense that it becomes the classical plane when $q$ is equal to 1. Let's begin with recalling some facts about the exterior algebra.

\subsection{$Z_3$-graded vector spaces}\label{sec2.1}

A $Z_3$-graded vector space is a vector space $V$ together with a decomposition $V = V_0\oplus V_1\oplus V_2$. Members of $V$ are called homogeneous elements. The {\it degree} (or {\it grade}) of a homogenous element $v\in V_m$ is denoted by $\varrho(v)=m$, $m\in Z_3=\{0,1,2 \, (\mbox{mod 3})\}$. An element in $V_m$ is of degree $m$.

\subsection{$Z_3$-graded algebras}\label{sec2.2}

A $k$-algebra $A$ is a $Z_3$-graded vector space $A$ together with multiplication such that $A_m\cdot A_l \subseteq A^{m+l}$ for all $m,l \in Z_3$.
Instead of linear maps, we use $k$-algebra maps for the definition of graded maps.

\subsection{$Z_3$-graded line $\tilde{\mathbb R}^{0|1}$}\label{sec2.3}

$Z_3$-graded Grassmann algebra with one generator over ${\mathbb C}$ denoted by ${\cal O}(\tilde{\mathbb R}^{0|1})$ is an associative algebra generated by the element $\theta$ such that
\begin{equation} \label{1}
\theta^3=0
\end{equation}
where the generator $\theta$ is of grade 1. Therefore, all elements of ${\cal O}(\tilde{\mathbb R}^{0|1})$ are polynomials of at most 2nd degree in $\theta$ and the basis of ${\cal O}(\tilde{\mathbb R}^{0|1})$ as a vector space is the set $\{\theta^i: \, i\in\{0,1,2\}\}$ where $\theta^0={\bf 1}$. Dimension of ${\cal O}(\tilde{\mathbb R}^{0|1})$ is thus $3^1=3$.

\subsection{$Z_3$-graded (exterior) plane $\tilde{\mathbb R}_q^{0|2}$}\label{sec2.4}

The exterior quantum plane (which is called the dual of quantum plane by Manin \cite{Manin1}) is generated by elements $\theta$ and $\varphi$ which satisfy the conditions
\begin{equation*}
\theta\otimes \varphi = - p \,\varphi\otimes \theta, \quad \theta\otimes \theta = 0 = \varphi\otimes \varphi
\end{equation*}
where $p$ is a nonzero complex number (instead of $q^{-1}$ we used $p$, to avoid confusion). Conventionally, the tensor product is not openly written as we did above. (Roughly, it is needed to put the wedge product between the coordinates of exterior plane, but it does not matter in the $Z_3$-graded case.)

When the plane is considered together with the exterior plane, we understand that the degree of coordinates of the plane is increased and they are nilpotent. Therefore, the new coordinates become anti-commutative. The resulting plane is called the $Z_2$-graded exterior plane.
So, it seems to be possible to construct a certain graded plane from a classical plane. In accordance with this opinion, a $Z_3$-graded exterior algebra can be constructed from a classical exterior algebra presented by

\noindent(a) increasing the power of nilpotency of its generators

\noindent(b) increasing degrees of some of its generators and

\noindent(c) imposing a $Z_3$-graded commutation relation on the generators.



Let $K\{\theta,\varphi\}$ be a free algebra where the coordinate $\theta$ is of degree 1 and the coordinate $\varphi$ is of degree 2. We also assume that $q$ is a cubic root of unity and $q\ne1$.

\begin{definition} 
Let $I_q$ be a two-sided ideal generated by the elements $\theta^3$, $\theta\cdot\varphi-q^2\theta\cdot\varphi$ and $\varphi^3$. The exterior plane  $\tilde{\mathbb R}_q^{0|2}$ with the function algebra
\begin{equation*}
{\cal O}(\tilde{\mathbb R}_q^{0|2}) = K\{\theta,\varphi\}/I_q
\end{equation*}
is called Z$_3$-graded exterior plane.
\end{definition}

By Definition 2.1, we have
\begin{equation} \label{2}
\tilde{\mathbb R}_q^{0|2} \ni \begin{pmatrix} \theta \\ \varphi \end{pmatrix}  \,\, \Longleftrightarrow \,\, \theta\cdot\varphi = q^2 \theta\cdot\varphi, \quad \theta^3=0=\varphi^3.
\end{equation}

Thus, all elements of ${\cal O}(\tilde{\mathbb R}_q^{0|2})$ are polynomials of at most 2nd degree in $\theta$ and $\varphi$ and the basis of
${\cal O}(\tilde{\mathbb R}_q^{0|2})$ is the set $\{\theta^i \cdot\varphi^j: \, i,j\in\{0,1,2\}\}$ where $\theta^0={\bf 1}=\varphi^0$. Dimension of
${\cal O}(\tilde{\mathbb R}_q^{0|2})$ is therefore $3^2=9$.

It is not hard to see the existence of representations that satisfy (\ref{2}); for instance, the two complex 3x3 matrix representations of the coordinate functions $\theta$ and $\varphi$
\begin{equation*}
\rho(\theta)=\begin{pmatrix} 0 & 1 & c_1 \\ 0 & 0 & q^2 \\ 0 & 0 & 0 \end{pmatrix}, \quad
\rho(\varphi)=\begin{pmatrix} 0 & q^2 & c_2 \\ 0 & 0 & 1 \\ 0 & 0 & 0 \end{pmatrix}
\end{equation*}
satisfy the relations (\ref{2}) for all numbers $c_1, c_2$.

\begin{definition}
Let ${\cal O}(\tilde{\mathbb R}_q^{*0|2})$ be the algebra with the generators $\xi$ and $x$ obeying the relations
\begin{equation} \label{3}
\xi \cdot x = x \cdot \xi
\end{equation}
where the coordinates $\xi$ and $x$ are of grade 2 and 0, respectively. We call ${\cal O}(\tilde{\mathbb R}_q^{*0|2})$ the algebra of functions on the Z$_3$-graded dual exterior plane $\tilde{\mathbb R}_q^{*0|2}$.
\end{definition}

\noindent{\bf Some conventions.}
\begin{enumerate}
\item Here, with the term "dual", we mean  to increase 1 (mod 3) the degree of an element. In fact, since $\varrho(\theta)=1$ and  $\varrho(\varphi)=2$ one has $\varrho(\mbox{dual of} \,\, \theta)=2$ (mod 3) and $\varrho(\mbox{dual of} \,\, \varphi)=0$ (mod 3).
\item As a conclusion, we modified or {\it deformed} the exterior plane which is freely generated by the two variables $\theta$, $\varphi$ subject to the trivial anti-commutation relations $\theta\varphi=-\varphi\theta$, $\theta^2=0=\varphi^2$. So, it can be called the $Z_3$-graded quantum exterior plane which is given in Definition 2.1.
\end{enumerate}

\section{$Z_3$-graded Hopf algebras}\label{sec3}

A unital algebra $H$ is called a Hopf algebra if there exist two algebra homomorphisms $\Delta:H\longrightarrow H\otimes H$, called the coproduct, and
$\epsilon:H\longrightarrow {\mathbb C}$, called the counit and an algebra anti-homomorphism $\kappa:H\longrightarrow H$, called the antipode of $H$, such that
\begin{eqnarray} 
&& (\Delta \otimes \textrm{id}) \circ \Delta = (\textrm{id} \otimes \Delta) \circ \Delta \label{4}\\
&& m \circ (\epsilon \otimes \mbox{id}) \circ \Delta = \mbox{id} = m \circ (\mbox{id} \otimes \epsilon) \circ \Delta   \label{5}\\
&& m \circ (\kappa \otimes \mbox{id}) \circ \Delta = \epsilon = m \circ (\mbox{id} \otimes \kappa) \circ \Delta \label{6}
\end{eqnarray}
where $m$ stands for multiplication in $H$.

A left quantum space (or briefly, a quantum space) for $H$ is an algebra ${\sf X}$ together with an algebra homomorphism
$\delta_L:{\sf X}\longrightarrow H\otimes{\sf X}$ such that
\begin{equation*}
(\mbox{id}\otimes\delta_L)\circ\delta_L=(\Delta\otimes\mbox{id})\circ\delta_L \quad \mbox{and} \quad (\epsilon\otimes\mbox{id})\circ\delta_L=\mbox{id}.
\end{equation*}
Also, $\delta_L$ is called the {\sf action} of $H$ on ${\sf X}$.

For these and more, we refer to the excellent books of Klimyk and Schm\"udgen \cite{Kli-sch} and Majid \cite{Majid1} for interested readers.

\subsection{The $Z_3$-graded Hopf algebra ${\cal O}(\tilde{\mathbb R}^{0|1})$}\label{sec3.1}

If $A$ and $B$ are two $Z_3$-graded algebra then the tensor product $A\otimes B$ is also a $Z_3$-graded algebra.
\begin{definition} 
For the $Z_3$-graded Grassmann algebra $A$, the product rule in the $Z_3$-graded Grassmann algebra $A\otimes A$ is defined by
\begin{equation} \label{7}
(a\otimes b)\cdot(c\otimes d) = \alpha(\varrho(b),\varrho(c)) \, a\cdot c\otimes b\cdot d
\end{equation}
where
\begin{equation} \label{8}
\alpha(\varrho(x),\varrho(y)) = \left\{\begin{array}{r@{,\quad}l} q^{\varrho(x)-\varrho(y)} & x\neq y \\[2ex] q^{\varrho(x)} & x=y \end{array}\right.
\end{equation}
and $q$ is a cubic root of unity.
\end{definition}

A $Z_3$-graded algebra is a $Z_3$-graded Hopf algebra if it equipped with co-maps given in (\ref{4})-(\ref{6}) such that the condition (\ref{7}) holds.

\begin{theorem}
The algebra ${\cal O}(\tilde{\mathbb R}^{0|1})$ is a $Z_3$-graded Hopf algebra with the definitions of coproduct, counit and antipode on the algebra ${\cal O}(\tilde{\mathbb R}^{0|1})$ as follows
\begin{equation} \label{9}
\Delta(\theta) = \theta \otimes {\bf 1} + {\bf 1}\otimes \theta, \quad \epsilon(\theta)=0, \quad \kappa(\theta) = - \theta.
\end{equation}
\end{theorem}

\noindent{\it Proof.} 
It is not hard to see that the properties (\ref{4})-(\ref{6}) are satisfied. Therefore, we only show that $\Delta(\theta^3)=0$ according to (\ref{1}). Using the rule given with (\ref{7}), we write
\begin{eqnarray*}
\Delta(\theta^2) &=& \theta^2 \otimes {\bf 1} - q^2\theta\otimes \theta + {\bf 1}\otimes \theta^2\\
\Delta(\theta^3) &=& \theta^3 \otimes {\bf 1} - q^2 \theta^2 \otimes \theta + \theta\otimes \theta^2 + q^2 \theta^2 \otimes \theta - \theta\otimes \theta^2 + {\bf 1}\otimes \theta^3
\end{eqnarray*}
which is the desired result with $\theta^3=0$. \hfill $\square$

\subsection{The $Z_3$-graded Hopf algebra ${\cal O}(\tilde{\mathbb R}_q^{0|2})$}\label{sec3.2}

The Theorem 3.1 can be extended to the case of two elements as follows:
\begin{theorem}
Define a coproduct $\Delta$, a counit $\epsilon$ for ${\cal O}(\tilde{\mathbb R}_q^{0|2})$ as
\begin{eqnarray} 
& \Delta(\theta) = \theta \otimes {\bf 1} + {\bf 1}\otimes \theta, & \Delta(\varphi) = \varphi \otimes {\bf 1} + {\bf 1}\otimes \varphi, \label{10}\\
& \epsilon(\theta)=0, & \epsilon(\varphi)=0. \label{11}
\end{eqnarray}
These maps supply ${\cal O}(\tilde{\mathbb R}_q^{0|2})$ with a $Z_3$-graded bialgebra structure.
\end{theorem}

\noindent{\it Proof.} 
It is easy to see that the properties (\ref{4}),(\ref{5}) are satisfied. Therefore, we only show that $\Delta(\theta\cdot\varphi - q^2 \varphi\cdot\theta)=0$, because it is shown that $\Delta(\theta^3)=0$ in Theorem 3.1. Using the rule given with (\ref{7}), we write
\begin{eqnarray*}
\Delta(\theta\cdot\varphi)
&=& \Delta(\theta)\cdot\Delta(\varphi) = \theta\cdot\varphi \otimes {\bf 1} + \theta\otimes \varphi + q^2 \varphi\otimes\theta + {\bf 1}\otimes \theta\cdot\varphi\\
\Delta(\varphi\cdot\theta)
&=& \Delta(\varphi)\cdot\Delta(\theta) = \varphi\cdot\theta \otimes {\bf 1} + \varphi\otimes\theta + q \theta\otimes \varphi + {\bf 1}\otimes \varphi\cdot\theta.
\end{eqnarray*}
So, according to the first relation in (\ref{2}) it must be $\Delta(\theta\cdot\varphi - q^2 \varphi\cdot\theta)=0$. \hfill $\square$

\begin{corollary}
An antipode $\kappa$ for ${\cal O}(\tilde{\mathbb R}_q^{0|2})$ is defined by
\begin{equation} \label{12}
\kappa(\theta) = - \theta, \quad \kappa(\varphi) = - \varphi.
\end{equation}
The antipode $\kappa$ of ${\cal O}(\tilde{\mathbb R}_q^{0|2})$ is an algebra anti-homomorphism and a coalgebra anti-homomorphism of
${\cal O}(\tilde{\mathbb R}_q^{0|2})$.
\end{corollary}

\noindent{\it Proof.} 
It is easy to check that the property (\ref{6}) holds. Let the linear map $m$ be multiplication on the algebra ${\cal O}(\tilde{\mathbb R}_q^{0|2})$ and the letter $\tau$ be the $Z_3$-graded flip operator given by
\begin{equation} \label{13}
\tau(a\otimes b) = q^{\varrho(a)\varrho(b)} b\otimes a.
\end{equation}
To prove the first claim, we must show that
\begin{equation} \label{14}
\kappa\circ m = m\circ \tau\circ (\kappa\otimes\kappa) \quad \mbox{and} \quad \kappa(1) = 1.
\end{equation}
An easy computation shows that
\begin{eqnarray*}
m\circ (\tau(\kappa\otimes\kappa)(\theta\otimes\varphi))
&=& m(\tau(\kappa(\theta)\otimes \kappa(\varphi))) = m(\tau(\theta\otimes \varphi)) \\
&=& m(q^2 \varphi\otimes \theta) = q^2 \varphi\cdot \theta = q^2 \kappa(\varphi)\cdot \kappa(\theta) \\
\kappa(m(\theta\otimes\varphi)) &=& \kappa(\theta\cdot\varphi).
\end{eqnarray*}
The antipode is calculated by the property
\begin{equation} \label{15}
\kappa(ab) = q^{\varrho(a)\varrho(b)} \, \kappa(b)\kappa(a)
\end{equation}
in terms of antipode of the generators. The second one follows at once from the identity   combined with the facts that
$m\circ (\kappa\otimes\mbox{id})\Delta({\bf 1})=\epsilon(1){\bf 1}$ and $\epsilon(1)=1$. To prove that the antipode is a coalgebra anti-homomorphism, we check the identities
\begin{equation} \label{16}
\Delta\circ\kappa = \tau\circ(\kappa\otimes\kappa)\circ\Delta \quad \mbox{and} \quad \epsilon\circ\kappa = \epsilon.
\end{equation}
Some computations give
\begin{eqnarray*}
\tau\circ(\kappa\otimes\kappa)\Delta(\theta)
&=& \tau(\kappa\otimes\kappa)(\theta \otimes {\bf 1} + {\bf 1}\otimes \theta) \\
&=& - \tau(\theta \otimes {\bf 1} + {\bf 1}\otimes \theta) \\
&=& - ({\bf 1}\otimes \theta + \theta \otimes {\bf 1}) \\
&=& \Delta(\kappa(\theta)).
\end{eqnarray*}
The other is trivial. \hfill $\square$

The formula (\ref{10}) indicates the action of the coproduct $\Delta$ only on the generators. The action of $\Delta$ on product of generators can be calculated by considering that $\Delta$ is an algebra homomorphism.

\section{Differential calculi on the Z$_3$-graded exterior plane}\label{sec4}

In this section, we will set up differential calculi on the $Z_3$-graded exterior plane $\tilde{\mathbb R}_q^{0|2}$. These calculi involve functions on this plane, differentials and differential forms.

One of these calculi which will be developed here is covariant with respect to the action of the $Z_3$-graded quantum group \cite{Celik3}.

\subsection{Differential algebra of an algebra} \label{sec4.1}

Differential calculi can be constructed on algebras. We begin with the definition of the $Z_3$-graded differential calculus. Let ${\cal A}$ be an arbitrary algebra with unity and $\Gamma^n({\cal A}):=\Gamma^n$ be a space of $n$-form and ${\cal A}$-bimodule.

\begin{definition}
A $Z_3$-graded differential calculus on the algebra ${\cal A}$ is a $Z_3$-graded algebra $\Gamma=\bigoplus_{n\ge0} \Gamma^n$ with a linear differential operator {\sf d} which defines the map ${\sf d}:{\cal A} \longrightarrow \Gamma$ of degree 1. A generalization of a usual differential calculus leads to the rules: for $f,g\in {\cal A}$.

\noindent(1) ${\sf d}(f\cdot g) = ({\sf d}f)\cdot g + q^{\varrho(f)} f\cdot ({\sf d}g)$,

\noindent(2) ${\sf d}^2(f\cdot g) = ({\sf d}^2f)\cdot g - q^{\varrho(f)-1} \, {\sf d}f \wedge {\sf d}g +
 q^{2\varrho(f)} \, f\cdot ({\sf d}^2g)$,

\noindent(3) ${\sf d}^3f = 0, \quad ({\sf d}^2f\ne0)$.
\end{definition}
The first and second conditions in Definition 4.1 are called the $Z_3$-graded Leibniz rules of degrees 1 and 2, respectively. A $Z_3$-graded differential algebra over ${\cal A}$ is a $Z_3$-graded algebra $\Gamma=\bigoplus_{n\ge0} \Gamma^n$, $\Gamma^0:={\cal A}$, with the linear map {\sf d} of degree 1 such that
${\sf d}^3=0$ and $Z_3$-graded Leibniz rules hold.

\vspace*{.2cm}
\noindent{\it Some conventions and assumptions}:
The Z$_3$-graded exterior algebra underlies a differential calculus on a smooth manifold with the differential {\sf d} satisfying ${\sf d}^3=0$. So, in order to construct a differential calculus on the Z$_3$-graded exterior algebra, a linear operator {\sf d} which acts on the functions of the Z$_3$-graded exterior algebra must be defined. For the definition, it is sufficient to define the action of {\sf d} on the components of the vector $\Theta$ in
$\tilde{\mathbb R}_q^{0|2}$ and on their products as follows: The operator {\sf d} is applied to $\theta$ produces a 1-form of degree 2, by definition. Similarly, application of {\sf d} to $\varphi$ produces  1-form whose $Z_3$-grade are 0. When the operator {\sf d} is applied to twice by iteration to $\theta$ and $\varphi$ it will produce a new entities which we will call a 2-form of degree 0 and 1, denoted by ${\sf d}^2\theta$ and ${\sf d}^2\varphi$, respectively. Finally, we require that ${\sf d}^3=0$.

\subsection{The structure of commutation relations} \label{sec4.2}

To obtain $q$-commutation relations between the generators of ${\cal O}(\tilde{\mathbb R}_q^{0|2})$ and their first order differentials, we combine elements $\theta$, $\varphi$ and their differentials ${\sf d}\theta$, ${\sf d}\varphi$ which are considered as elements of a space
$\Gamma^1({\cal O}(\tilde{\mathbb R}_q^{0|2})):= \Gamma^1$ of 1-forms. Let us allow a multiplication of the differentials by the elements of ${\cal O}(\tilde{\mathbb R}_q^{0|2})$ from the left and from the right so that by the definition of the multiplications the resulting 1-form belongs to $\Gamma^1$ again. This means that $\Gamma^1$ is an ${\cal O}(\tilde{\mathbb R}_q^{0|2})$-bimodule.

In order to obtain the $Z_3$-graded commutation relations between the elements of the algebra
${\cal O}(\tilde{\mathbb R}_q^{0|2})$ and their first and second order differentials, we will use partially the consistency of calculus \cite{WZ}. Consistency of the calculus does not allow making direct proof of this proposition. Because, second order differentials appear when the differential operator {\sf d} applied to the first order differentials. So, it is necessary to add a set of second order differentials to the first order differentials. Appearance of higher order differentials is a peculiar property of the $Z_3$-graded calculus. The second order differentials
${\sf d}^2\theta$ and ${\sf d}^2\varphi$ are considered as generators of a space
$\Gamma^2({\cal O}(\tilde{\mathbb R}_q^{0|2}))\doteq \Gamma^2$ of 2-forms.

\begin{lemma} 
The generators of ${\cal O}(\tilde{\mathbb R}_q^{0|2})$ with the first order differentials satisfy the following $Z_3$-graded commutation relations
\begin{eqnarray} \label{17}
& \theta\cdot{\sf d}\theta = Q_1 \, {\sf d}\theta\cdot\theta,
& \theta\cdot{\sf d}\varphi = C_1 \, {\sf d}\varphi\cdot\theta + C_2 \, {\sf d}\theta\cdot\varphi, \nonumber \\
& \varphi\cdot{\sf d}\varphi = Q_2 \, {\sf d}\varphi \cdot \varphi,
& \varphi\cdot{\sf d}\theta = (C_2+q^2) \, {\sf d}\theta\cdot\varphi + (C_1-q) \, {\sf d}\varphi\cdot\theta,
\end{eqnarray}
where $C_1$ and $C_2$ are complex constants and $Q_1\in\{1,q\}$ and $Q_2\in\{1,q^2\}$ for each case.
\end{lemma}

\noindent{\it Proof.} 
In general the coordinates will not commute with their differentials. Here we consider the case when the commutation relations of the coordinates with the differentials are bilinear. (Contrary situation will be give in Theorem 4.5) Therefore, one assumes that the possible commutation relations of the generators with their first order differentials are of the form
\begin{eqnarray} \label{18}
\theta \cdot {\sf d}\theta &=& Q_1 \, {\sf d}\theta \cdot \theta, \quad \theta \cdot {\sf d}\varphi = C_1 \, {\sf d}\varphi \cdot \theta + C_2 \, {\sf d}\theta \cdot \varphi, \nonumber \\
\varphi \cdot {\sf d}\varphi &=& Q_2 \, {\sf d}\varphi \cdot \varphi, \quad \varphi \cdot {\sf d}\theta = C_3 \, {\sf d}\theta \cdot \varphi + C_4 \, {\sf d}\varphi \cdot \theta
\end{eqnarray}
where the constants $Q_i$ and $C_j$ are possible related to $q$.

Since $\theta^3=0$, if we apply the differential operator {\sf d} on both sides of this equality we get
\begin{equation*}
0 = {\sf d}(\theta^3) = {\sf d}\theta\cdot \theta^2 + q \theta\cdot {\sf d}(\theta^2) = {\sf d}\theta\cdot \theta^2 + q \theta\cdot ({\sf d}\theta\cdot \theta + q\theta\cdot {\sf d}\theta).
\end{equation*}
So one obtains $1+qQ_1+q^2Q_1^2=0$ by which we can choose either $Q_1=1$ or $Q_1=q$. Similarly we find $Q_2=1$ or $Q_2=q^2$. If we now apply the operator {\sf d} from left to first relation in (\ref{2}), one has
\begin{equation} \label{19}
C_3=C_2+q^2 \quad \mbox{and} \quad C_4=C_1-q
\end{equation}
which completes the proof. \hfill $\square$

\noindent{\bf Note.} If we multiply the equation $\theta\cdot\varphi - q^2 \theta\cdot\varphi=0$ from right with
${\sf d}\theta$ and ${\sf d}\varphi$ and use the second and fourth relations in (\ref{18}) we obtain $C_2C_4=0$, $(C_1-q^2Q_1)C_4=0$ and $(Q_2-q^2C_3)C_2=0$. So, the relations (\ref{17}) have three forms with respect to $C_1$ and $C_2$.

We now assume that commutation relation between the first order differentials of the generators of $O(\tilde{\mathbb R}_q^{0|2})$ is of the form
\begin{equation} \label{20}
{\sf d}\theta\wedge{\sf d}\varphi = Q_3 \, {\sf d}\varphi\wedge{\sf d}\theta
\end{equation}
where the constant $Q_3$ is possible related to $q$.

If we apply the differential {\sf d} from left to the relations (\ref{18}) with respect to the $Z_3$-graded Leibniz
rule we get, for example,
\begin{equation} \label{21}
\theta\cdot{\sf d}^2\varphi = q^2C_1{\sf d}^2\varphi \cdot \theta + q^2C_2{\sf d}^2\theta\cdot\varphi + K_1{\sf d}\theta\wedge {\sf d}\varphi,
\end{equation}
where $K_1=q^2(C_1Q_3^{-1} + q^2C_2-1)$.

\noindent{\bf Remark 1.} The relation (\ref{21}) (and others) is not homogeneous in the sense that the commutation relations of the coordinates with second order differentials include first order differentials as well.

To make them homogenous, we need the coefficients of the first order terms to be zero. Then, by applying the operator {\sf d} to the relations (\ref{20}) and (\ref{21}) from left, we have
\begin{equation*}
Q_3=1 \quad \mbox{and} \quad C_2=q(1-C_1).
\end{equation*}
In other hand, the equation $0=(\theta\cdot\varphi-q^2\theta\cdot\varphi){\sf d}^2\theta$ gives
\begin{equation} \label{22}
(1-C_1)(C_1-q)=0.
\end{equation}

For $C_1=q$, we have

\begin{theorem} 
Commutation relations of the generators of ${\cal O}(\tilde{\mathbb R}_q^{0|2})$ and their first and second order differentials are of the form
\begin{eqnarray} 
&\theta\cdot{\sf d}\theta = q \, {\sf d}\theta\cdot\theta, &
 \varphi\cdot{\sf d}\theta = q \, {\sf d}\theta\cdot\varphi, \nonumber \\
&\theta\cdot{\sf d}\varphi = q \, {\sf d}\varphi\cdot\theta + (q-q^2) \, {\sf d}\theta\cdot\varphi, &
 \varphi\cdot{\sf d}\varphi = {\sf d}\varphi\cdot\varphi, \label{23}\\
&\theta\cdot{\sf d}^2\theta = {\sf d}^2\theta\cdot\theta, &
 \varphi\cdot{\sf d}^2\theta = q^2 \, {\sf d}^2\theta\cdot\varphi,\nonumber \\
&\theta\cdot{\sf d}^2\varphi = {\sf d}^2\varphi \cdot \theta + (1-q) \, {\sf d}^2\theta\cdot\varphi, &
 \varphi\cdot{\sf d}^2\varphi = q \, {\sf d}^2\varphi\cdot\varphi.\label{24}
\end{eqnarray}
\end{theorem}

\noindent{\bf Remark 2.} The relations (\ref{23}) (and thus (\ref{24})) are not invariant under $Z_3$-graded quantum group relations that will be defined in the next section. On the other hand, the relations obtained by taking $C_1=1$ will be invariant (see, Theorem 4.5).

The relations in the following theorem can be obtained with twice applying the differential {\sf d} to the second or third relation in (\ref{24}).
\begin{theorem} 
The commutation relation between the second order differentials is as follows
\begin{equation} \label{25}
{\sf d}^2\theta\wedge {\sf d}^2\varphi = q \, {\sf d}^2\varphi\wedge {\sf d}^2\theta.
\end{equation}
\end{theorem}

\noindent{\bf Remark 3.} It is not difficult to see that, the cubes of ${\sf d}\theta$ and ${\sf d}\varphi$ are both central elements of the whole differential algebra.

\noindent{\bf Remark 4.} It can be easily shown that
\begin{equation*}
{\sf d}({\sf d}\theta\wedge{\sf d}\theta\wedge{\sf d}\theta) = 0 = {\sf d}({\sf d}\varphi\wedge{\sf d}\varphi\wedge{\sf d}\varphi).
\end{equation*}
However, we will assume that, to be consistent,
\begin{equation} \label{26}
{\sf d}\theta\wedge{\sf d}\theta\wedge{\sf d}\theta \ne 0 \quad \mbox{and} \quad {\sf d}\varphi\wedge{\sf d}\varphi\wedge{\sf d}\varphi \ne 0.
\end{equation}
As it will be observed in Theorem 5.2, the cubes of partial derivatives cannot be zero unless
$({\sf d}\theta)^3: = {\sf d}\theta\wedge{\sf d}\theta\wedge{\sf d}\theta \ne 0$ and $({\sf d}\varphi)^3 \ne 0$.

\noindent{\bf Remark 5.} If we introduce the first order differentials of the generators of
${\cal O}(\tilde{\mathbb R}_q^{0|2})$ as
\begin{equation} \label{27}
\xi = {\sf d}\theta, \quad x = {\sf d}\varphi
\end{equation}
the relation (\ref{20}) coincides with (\ref{3}) where $Q_3$ is equal to 1.

Consequently, we set up $Z_3$-graded differential schemas with the relations satisfied by the elements of the set $\{\theta,\varphi,{\sf d}\theta,{\sf d}\varphi,{\sf d}^2\theta,{\sf d}^2\varphi\}$.

\subsection{Covariance of the calculus}\label{sec4.4}

We begin with the definition of the $Z_3$-graded quantum group $\widetilde{GL}_q(2)$. A detail discussion is given in the another paper \cite{Celik3}. Let $a$, $\beta$, $\gamma$, $d$ be elements of the algebra ${\cal O}(\tilde{M}(2))$. We also assume that the generators $a$ and $d$ are of degree 0, the generators $\gamma$ and $\beta$ are of degree 1 and 2, respectively. Then we have

\begin{definition} 
The $Z_3$-graded algebra ${\cal O}(\tilde{M}_q(2))$ is the quotient of the free algebra $k\{a,\beta,\gamma,d\}$ by the two-sided ideal $J_q$ generated by
the six relations
\begin{eqnarray} 
a \beta &=& \beta a, \quad \beta \gamma = \gamma \beta, \quad d\beta = \beta d, \label{28}\\
a \gamma &=& q \gamma a, \quad d\gamma = q^2 \gamma d, \label{29}\\
ad &=& da + (q-1) \beta \gamma, \label{30}
\end{eqnarray}
where $q$ is a cubic root of unity.
\end{definition}
By relation (\ref{30}), we have
\begin{equation} \label{31}
D_q := ad - q \beta \gamma = da - \beta \gamma.
\end{equation}
This element of ${\cal O}(\tilde{M}_q(2))$ is called the {\it $Z_3$-graded quantum determinant}.
Using the determinant $D_q$ belonging to  the algebra ${\cal O}(\tilde{M}_q(2))$, we can define a {\it new} Hopf algebra adding an inverse $t^{-1}$ to ${\cal O}(\tilde{M}_q(2))$. Let ${\cal O}(\widetilde{GL}_q(2))$ be the quotient of the algebra ${\cal O}(\tilde{M}_q(2))$ by the two-sided ideal generated by the element $tD_q-1$.

The $Z_3$-graded algebra ${\cal O}(\tilde{\mathbb R}_q^{0|2})$ is a (left) quantum space for the $Z_3$-graded quantum group $\widetilde{GL}_q(2)$, where its action $\delta_L$ is given by
\begin{equation} \label{32}
\delta_L(\theta)=a\otimes\theta+\beta\otimes\varphi, \quad \delta_L(\varphi)=\gamma\otimes\theta + d\otimes\varphi.
\end{equation}
More precisely, $\delta_L$ has to be extended to an algebra homomorphism of ${\cal O}(\tilde{\mathbb R}_q^{0|2})$ into
${\cal O}(\widetilde{GL}_q(2))\otimes {\cal O}(\tilde{\mathbb R}_q^{0|2})$.

\begin{definition} 
If ${\cal A}$ is a $Z_3$-graded algebra, then the product rule in the $Z_3$-graded algebra ${\cal A}\otimes {\cal A}$ is defined by
\begin{equation} \label{33}
(a_1\otimes a_2)(a_3\otimes a_4) = q^{\varrho(a_2)\varrho(a_3)} a_1a_3\otimes a_2a_4
\end{equation}
where $a_i$'s are homogeneous elements in the algebra ${\cal A}$.
\end{definition}

\begin{lemma} 
The linear mapping $\delta_L: {\cal O}(\tilde{\mathbb R}_q^{0|2})\longrightarrow {\cal O}(\widetilde{GL}_q(2))\otimes {\cal O}(\tilde{\mathbb R}_q^{0|2})$ is well-defined.
\end{lemma}

\noindent{\it Proof.} 
For this, we check that $\delta_L$ leaves invariant the relations $\theta\varphi=q^2\varphi\theta$ and $\theta^3=0=\varphi^3$. Since
\begin{eqnarray*}
\delta_L(\theta\varphi)
&=& qa\gamma\otimes\theta^2 + ad\otimes \theta\varphi + q^2\beta\gamma\otimes\varphi\theta + \beta d\otimes\varphi^2 \\
\delta_L(\varphi\theta)
&=& \gamma a\otimes\theta^2 + da\otimes\varphi\theta + q^2\gamma\beta\otimes\theta\varphi + qd\beta\otimes\varphi^2
\end{eqnarray*}
we have
\begin{eqnarray*}
\delta_L(\theta\varphi)-q^2\delta_L(\varphi\theta)
&=& q(a\gamma - q\gamma a)\otimes\theta^2 + (ad - q\gamma\beta)\otimes\theta\varphi \\
& & + q^2(\beta\gamma - da)\otimes\varphi\theta + (\beta d - d\beta)\otimes\varphi^2 \\
&=& 0
\end{eqnarray*}
as expected. Similarly, it is not difficult to check that $\delta_L(\theta^3)=0=\delta_L(\varphi^3)$. \hfill $\square$

In Hopf algebra terminology, the first suggestion of Definition 4.3 yields the following.
\begin{theorem}[\cite{Celik3}] 
The algebra ${\cal O}(\tilde{\mathbb R}_q^{0|2})$ is a left comodule algebra of the bialgebra ${\cal O}(\widetilde{GL}_q(2))$ with left coaction $\delta_L$ given by (\ref{32}).
\end{theorem}

\begin{definition} 
Let $A$ be a Hopf algebra, $V$ be a quantum space for $A$ with action $\delta_L$. A $Z_3$-graded differential calculus $(\Gamma,{\sf d})$ over $V$ is said to be {\sf left-covariant} with respect to $A$ if there exists an algebra homomorphism  $\Delta_L: \Gamma\longrightarrow A\otimes \Gamma$ which is a left coaction of $A$ on $\Gamma$  such that

\noindent(1) $\Delta_L|_V=\delta_L$ (that is, $\Delta_L(u)=\delta_L(u)$ for $u\in V$)

\noindent(2) $\Delta_L\circ{\sf d}=(\varrho\otimes{\sf d})\circ\Delta_L$ and
$\Delta_L\circ{\sf d}^2=(\varrho^2\otimes{\sf d}^2)\circ\Delta_L$.
\end{definition}

Let $\{{\sf d}\theta, {\sf d}\varphi, {\sf d}^2\theta, {\sf d}^2\varphi\}$ be a free right ${\cal O}(\tilde{\mathbb R}_q^{0|2})$-module basis of $\Gamma$. The proof of the following theorem follows from the Definition 4.4.

\begin{theorem} 
There exists a differential calculus $(\Gamma,{\sf d})$ on the algebra ${\cal O}(\tilde{\mathbb R}_q^{0|2})$ which are left-covariant
with respect to the $Z_3$-graded Hopf algebra ${\cal O}(\widetilde{GL}_q(2))$. The ${\cal O}(\tilde{\mathbb R}_q^{0|2})$-bimodule
structure of this calculus is described by the relations
\begin{eqnarray} 
&\theta\cdot{\sf d}\theta = q \, {\sf d}\theta\cdot\theta, &
 \varphi\cdot{\sf d}\theta = q^2 \, {\sf d}\theta\cdot\varphi +(1-q) \, {\sf d}\varphi\cdot\theta, \nonumber \\
&\theta\cdot{\sf d}\varphi = {\sf d}\varphi\cdot\theta, &
 \varphi\cdot{\sf d}\varphi = {\sf d}\varphi\cdot\varphi, \label{34}\\
&\theta\cdot{\sf d}^2\theta = {\sf d}^2\theta\cdot\theta, &
  \varphi\cdot{\sf d}^2\theta = {\sf d}^2\theta\cdot\varphi + (q-q^2) \, {\sf d}^2\varphi\cdot\theta, \nonumber\\
&\theta\cdot{\sf d}^2\varphi = q^2 \, {\sf d}^2\varphi\cdot\theta, \quad &
 \varphi\cdot{\sf d}^2\varphi = q \, {\sf d}^2\varphi\cdot\varphi, \label{35}\\
&{\sf d}\theta\wedge{\sf d}\varphi = {\sf d}\varphi\wedge{\sf d}\theta,
&{\sf d}^2\theta\wedge{\sf d}^2\varphi = q \, {\sf d}^2\varphi\wedge{\sf d}^2\theta.\label{36}
\end{eqnarray}
\end{theorem}

The above theorem has been obtained with the demand that the only linear relation between $\theta\cdot {\sf d}\theta$, $\theta\cdot {\sf d}\varphi$,
$\varphi\cdot {\sf d}\theta$, $\varphi\cdot {\sf d}\varphi$ and ${\sf d}\theta\cdot \theta$, ${\sf d}\varphi\cdot \theta$, ${\sf d}\theta\cdot \varphi$,
${\sf d}\varphi\cdot \varphi$. However, even if the demand is not considered, there exists a left covariant differential calculus:

\begin{theorem} 
There exists a differential calculus $(\Gamma,{\sf d})$ on the algebra ${\cal O}(\tilde{\mathbb R}_q^{0|2})$ which are left-covariant
with respect to the $Z_3$-graded Hopf algebra ${\cal O}(\widetilde{GL}_q(2))$. The ${\cal O}(\tilde{\mathbb R}_q^{0|2})$-bimodule
structure of this calculus is described by the relations
\begin{eqnarray*} 
\theta\cdot{\sf d}\theta &=& {\sf d}\theta\cdot\theta, \quad
 \varphi\cdot{\sf d}\theta = q \, {\sf d}\theta\cdot\varphi +(q^2-1) \, {\sf d}\varphi\cdot\theta, \nonumber \\
\theta\cdot{\sf d}\varphi &=& q^2 \, {\sf d}\varphi\cdot\theta, \quad
 \varphi\cdot{\sf d}\varphi = q^2 \, {\sf d}\varphi\cdot\varphi, \\
\theta\cdot{\sf d}^2\theta &=& q^2 \, {\sf d}^2\theta\cdot\theta + (q-q^2) \, {\sf d}\theta\wedge{\sf d}\theta, \\
\theta\cdot{\sf d}^2\varphi &=& q \, {\sf d}^2\varphi\cdot\theta + (q-q^2) \, {\sf d}\varphi\wedge{\sf d}\theta, \\
\varphi\cdot{\sf d}^2\theta &=& q^2 \, {\sf d}^2\theta\cdot\varphi + (1-q) \, {\sf d}^2\varphi\cdot\theta + (1-q) \, {\sf d}\theta\wedge{\sf d}\varphi, \\
\varphi\cdot{\sf d}^2\varphi &=& q \, {\sf d}^2\varphi\cdot\varphi + (1-q) \, {\sf d}\varphi\wedge{\sf d}\varphi, \\
{\sf d}\theta\wedge{\sf d}\varphi &=& {\sf d}\varphi\wedge{\sf d}\theta, \quad {\sf d}^2\theta\wedge{\sf d}^2\varphi = q \, {\sf d}^2\varphi\wedge{\sf d}^2\theta.
\end{eqnarray*}
\end{theorem}

\section{Relations with the partial derivatives}\label{sec5}

In this section, we will continue with the relations given in Theorem 4.4.
To obtain the commutation relations between the generators of ${\cal O}(\tilde{\mathbb R}_q^{0|2})$ and derivatives, let us introduce the partial derivatives of the generators of the algebra. Since $\Gamma$ is a left covariant differential calculus, for any element $\alpha$ in ${\cal O}(\tilde{\mathbb R}_q^{0|2})$ there are uniquely determined elements $\partial_k(\alpha)\in {\cal O}(\tilde{\mathbb R}_q^{0|2})$ such that
\begin{equation*}
{\sf d}\alpha = \sum_k {\sf d}\alpha_k \partial_k(\alpha), \quad \partial_k:=\partial/\partial\alpha_k
\end{equation*}
where $\partial_k\alpha_l=\delta_{kl}$. The mappings $\partial_k:{\cal O}(\tilde{\mathbb R}_q^{0|2})\longrightarrow {\cal O}(\tilde{\mathbb R}_q^{0|2})$ are called the {\it partial derivatives} of the calculus $\Gamma$.


\subsection{Relations between the coordinates and the partial \\ derivatives}\label{sec5.1}

The following theorem gives the relations between the elements of ${\cal O}(\tilde{\mathbb R}_q^{0|2})$ and their partial derivatives.

\begin{theorem} 
The relations of the coordinates with their partial derivatives are as follows
\begin{eqnarray} \label{37}
\partial_\theta \, \theta &=& 1+q^2 \, \theta \, \partial_\theta, \quad \partial_\theta \, \varphi = q \, \varphi \, \partial_\theta, \nonumber \\
\partial_\varphi \, \theta &=& q \, \theta \, \partial_\varphi, \quad
\partial_\varphi \, \varphi = 1 + q^2 \, \varphi \, \partial_\varphi + (q^2-1) \, \theta \, \partial_\theta.
\end{eqnarray}
\end{theorem}

\noindent{\it Proof.} 
We know that the differential {\sf d} can be written in terms of the differentials and partial derivatives as
\begin{eqnarray} \label{38}
{\sf d}f = ({\sf d}\theta \, \partial_\theta + {\sf d}\varphi \, \partial_\varphi)f
\end{eqnarray}
where $f$ is a differentiable function. For consistency, the degrees of the derivatives $\partial_\theta$ and $\partial_\varphi$ should be 2, 1, respectively. Now, if we replace $f$ with $\theta f$ on the left hand side of the equality in (\ref{38}), we get
\begin{eqnarray*}
{\sf d}(\theta f) = {\sf d}\theta \, f+ q \, \theta \, {\sf d}f = \left[{\sf d}\theta \, (1+q^2 \, \theta \, \partial_\theta) + {\sf d}\varphi \, (q \, \theta \, \partial_\varphi)\right]f. \nonumber
\end{eqnarray*}
On the other hand, the right hand side of the equality in (\ref{38}) is of the form
\begin{eqnarray*}
{\sf d}(\theta f) = [{\sf d}\theta \, (\partial_\theta \, \theta) + {\sf d}\varphi \, (\partial_\varphi \, \theta)]f.
\end{eqnarray*}
Now, by comparing the right hand sides of these two equalities according to the differentials we obtain some of the relations in (\ref{37}). Other relations can be found similarly. \hfill $\square$

\subsection{Relations between the partial derivatives}\label{sec5.2}

Since we have the condition(s) ${\sf d}^3=0$ (and ${\sf d}^2\ne0$) in $Z_3$-graded calculus, it requires some tedious computations to find the relations of partial derivatives. Once the operator {\sf d} is applied to the equation (\ref{38}), one obtains
\begin{eqnarray*}
{\sf d}^2f
&=& \left[{\sf d}^2\theta \,\partial_\theta + {\sf d}^2\varphi \,\partial_\varphi + ({\sf d}\theta \wedge {\sf d}\varphi)(\partial_\theta \partial_\varphi +
     q^2\partial_\varphi\partial_\theta)\right. \\& & \left. + q^2({\sf d}\theta \wedge {\sf d}\theta)\partial_\theta^2 + ({\sf d}\varphi \wedge {\sf d}\varphi) \,\partial_\varphi^2\right]f.
\end{eqnarray*}
Next, the operator {\sf d} will be re-applied to both sides of this equation and then the conditions in (\ref{26}) will be used to obtain the desired relations. As a result, we have the following theorem:


\begin{theorem} 
The relations between partial derivatives are of the form
\begin{eqnarray} \label{39}
\partial_\varphi \partial_\theta = q^2 \, \partial_\theta \partial_\varphi, \quad \partial_\theta^3 = 0 = \partial_\varphi^3.
\end{eqnarray}
\end{theorem}
The last two equalities in (\ref{39}) are consequences of assumptions given in (\ref{26}).

\subsection{Relations of differentials with the partial derivatives}\label{sec5.3}

To complete the schema, we need the relations that partial derivatives are satisfied with first and second order differentials.

\begin{theorem} 
The relations of partial derivatives with first and second order differentials are as follows
\begin{eqnarray} \label{40}
{\sf d}\theta \, \partial_\theta &=& q \, \partial_\theta \, {\sf d}\theta + (q^2-1) \, \partial_\varphi \,  {\sf d}\varphi, \quad
 {\sf d}\varphi \, \partial_\theta = \partial_\theta \, {\sf d}\varphi, \nonumber \\
{\sf d}\theta \, \partial_\varphi &=& q^2 \, \partial_\varphi \, {\sf d}\theta, \quad {\sf d}\varphi \, \partial_\varphi = \partial_\varphi \, {\sf d}\varphi,
\end{eqnarray}
and
\begin{eqnarray} \label{41}
{\sf d}^2\theta \, \partial_\theta &=& \partial_\theta \, {\sf d}^2\theta + (1-q) \, {\sf d}^2\varphi \, \partial_\varphi, \quad
 {\sf d}^2\varphi \,  \partial_\theta = q^2 \, \partial_\theta \, {\sf d}^2\varphi, \nonumber\\
{\sf d}^2\theta \, \partial_\varphi &=& \partial_\varphi \, {\sf d}^2\theta, \quad {\sf d}^2\varphi \, \partial_\varphi = q \, \partial_\varphi \,  {\sf d}^2\varphi.
\end{eqnarray}
\end{theorem}

\noindent{\it Proof.} 
To obtain the relations in (\ref{40}), we assume that the following relations are satisfied with partial derivatives and first order differentials:
\begin{eqnarray*}
\partial_\theta \, {\sf d}\theta &=& A_1 \, {\sf d}\theta \, \partial_\theta + A_2 \, {\sf d}\varphi \, \partial_\varphi, \quad
\partial_\theta \, {\sf d}\varphi = A_3 \, {\sf d}\varphi \, \partial_\theta + A_4 \, {\sf d}\theta \, \partial_\varphi, \\
\partial_\varphi \, {\sf d}\varphi &=& A_5 \, {\sf d}\varphi \, \partial_\varphi + A_6 \, {\sf d}\theta \, \partial_\theta, \quad
\partial_\varphi \, {\sf d}\theta = A_7 \, {\sf d}\theta \, \partial_\varphi + A_8 \, {\sf d}\varphi \, \partial_\theta.
\end{eqnarray*}
The constants $A_i$'s will be determined in terms of the parameter $q$. To find them, we apply the operator  $\partial_\theta$ to the relations (\ref{34}). Using the fact
\begin{eqnarray*}
\partial_i (X^j {\sf d}X^k) = \delta^i_j \,{\sf d}X^k
\end{eqnarray*}
where $\partial_1 = \partial_\theta$, $\partial_2 = \partial_\varphi$, $X^1 = \theta$ and $X^2=\varphi$, after from some calculations we find $A_1=q^2$, $A_2=q^2-q$, $A_3=1$, $A_4=0$, $A_5=q$, $A_6=0$, $A_7=1$ and $A_8=0$. Other relations in (\ref{40}) can be determined in a similar manner. The relations in (\ref{41}) can also be obtained by using the same idea. \hfill $\square$

\noindent Note that the relations of the differential {\sf d} with $\partial_\theta$ and $\partial_\varphi$ are
\begin{eqnarray} \label{42}
\partial_k \, {\sf d} = q^{1+2\varrho(\partial_k)} \, {\sf d} \, \partial_k \quad \mbox{and} \quad \partial_k \, {\sf d}^2 = q^{2+\varrho(\partial_k)} \, {\sf d}^2 \, \partial_k.
\end{eqnarray}

\section{$Z_3$-graded Weyl and Heisenberg algebras}\label{sec6}

An appropriate interaction between the coordinates and partial derivatives can be defined as in \cite{WZ} and
also \cite{Giaq}.

\begin{definition} 
The $Z_3$-graded Grassmann-Weyl algebra $\tilde{\cal A}_q(2)$ is the unital algebra generated by $\theta$, $\varphi$ and $\partial_\theta$, $\partial_\varphi$ which satisfy the relations (\ref{2}), (\ref{37}) and (\ref{39}).
\end{definition}
Obviously, the algebra $\tilde{\cal A}_q(2)$ is a PI ring and is not primitive.

The following conclusion asserts that relations in (\ref{39}) give no new relations between the $\theta$, $\varphi$ and $\partial_\theta$, $\partial_\varphi$ other than those given in (\ref{2}) and (\ref{37}).

\begin{corollary} 
We have that, with $I_q$ given in Definition 2.1, $\partial_uI_q = I_q\partial_u$, for $u\in \{\theta,\varphi\}$.
\end{corollary}

\begin{definition} 
Let $\cal A$ be a $Z_3$-graded algebra. A conjugate-linear map $\xi \mapsto \xi^\star$ of degree zero is called a $Z_3$-graded involution on $\cal A$ if
\begin{equation} \label{43}
(\xi_1\xi_2)^\star = q^{\varrho(\xi_1)\varrho(\xi_2)} \xi_2^\star \xi_1^\star, \quad (\xi^\star)^\star = \xi
\end{equation}
for all $\xi, \xi_1, \xi_2 \in A$. The pair $(\cal A,\star)$ is called a $Z_3$-graded $\star$-algebra.
\end{definition}

We now would like to define a $Z_3$-graded Heisenberg algebra. Therefore, let us choose the hermitean conjugations of the coordinates $\theta$, $\varphi$ and the derivatives $\partial_\theta$ and $\partial_\varphi$ as $\theta^\star = \theta$, $\varphi^\star = \varphi$ and $\partial_\theta^\star = -\partial_\theta$, $\partial_\varphi^\star = -q^2\partial_\varphi$. Then, the operation $\star$ leaves the relations (\ref{2}), (\ref{37}) and (\ref{39}) as invariant. This involution allows us to define the hermitean operators
\begin{equation} \label{44}
\hat{\theta} = \theta, \quad \hat{\varphi} = \varphi, \quad \hat{p}_\theta = i \partial_\theta, \quad \hat{p}_\varphi = i q \partial_\varphi.
\end{equation}
Let us denote the new algebra by $\tilde{\cal H}_q(2)$ whose generators belong to the set $\{\hat{p}_\theta,\hat{p}_\varphi,\hat{\theta},\hat{\varphi}\}$.
Then, the following proposition can be proven directly.
\begin{theorem} 
The $Z_3$-graded algebra $\tilde{\cal H}_q(2)$ is a $Z_3$-graded Lie algebra on \\
$span_{\mathbb C}\{\hat{p}_\theta,\hat{p}_\varphi,\hat{\theta},\hat{\varphi},i{\bf 1}\}$ with the following commutation relations
\begin{eqnarray} \label{45}
\hat{\theta}^3 = 0, & \hat{\varphi}^3 = 0, & [\hat{\theta},\hat{\varphi}]_{q^2} = 0, \nonumber\\
\hat{p}_\theta^3=0, & [\hat{p}_\theta, \hat{\theta}]_{q^2} = i {\bf 1}, & [\hat{p}_\theta, \hat{\varphi}]_q = 0, \\
\hat{p}_\varphi^3=0, & [\hat{p}_\varphi,\hat{\theta}]_q=0, & [\hat{p}_\varphi,\hat{\varphi}]_q = iq {\bf 1} + (1-q) \hat{\theta}\hat{p}_\theta. \nonumber
\end{eqnarray}
\end{theorem}

A $Z_3$-graded derivation operator acting on the algebra ${\cal O}(\tilde{\mathbb R}_q^{0|2})$ is defined by
\begin{equation} \label{46}
D = \theta\partial_\theta + \varphi\partial_\varphi.
\end{equation}
As an element of the $Z_3$-graded Weyl algebra, the operator $D$ is the $Z_3$-graded Euler derivation of grade 0.
An easy computation shows that $D + D^\star = {\bf 1}$ and $D$ is a normal operator, but not invertible. However, the operator $E:={\bf 1}+(q^2-1) D$ is invertible and by (\ref{37}),(\ref{39}), it is easy to check that $E\theta_i = q^2\theta_i E$ and
$E \partial_{\theta_i} = q \partial_{\theta_i} E$ and $E^\star = q E$ so that $E$ is a normal operator.

\section{An $R$-Matrix Formalism}\label{sec7}

Using the relations in (\ref{34}), we can find an $R$-matrix satisfying the $Z_3$-graded Yang-Baxter equation. We assume that an $R$-matrix is associated with the Z$_3$-graded plane $\tilde{\mathbb R}_q^{0|2}$. Then, we can express the commutation relations, (\ref{34}), between the coordinates and their differentials in the form
\begin{eqnarray} \label{47}
q^{\varrho(\theta_i)} \theta_i {\sf d}\theta_j = q \sum_{k,l=1}^2 \hat{R}_{ij}^{kl} \, {\sf d}\theta_k \theta_l,
\end{eqnarray}
where $\hat{R}=\underline{P}R$ with the $Z_3$-graded permutation matrix $\underline{P}$. If we compare it with the relations (\ref{34}), we see that, except the zero entries, $\hat{R}^{11}_{11}=q=\hat{R}^{22}_{22}$, $\hat{R}^{12}_{21}=1=\hat{R}^{21}_{12}$ and $\hat{R}^{21}_{21}=q-q^2$. So, $Z_3$-graded $\hat{R}$-matrix takes the following form
\begin{eqnarray} \label{48}
\hat{R} = \begin{pmatrix}
q & 0 & 0 & 0 \\
0 & 0 & 1 & 0 \\
0 & 1 & q-q^2 & 0 \\
0 & 0 & 0 & q
\end{pmatrix}=(\hat{R}_{ij}^{kl}).
\end{eqnarray}

The commutation rules of the generators, $\Theta=(\theta_i)$, of function algebra on the plane $\tilde{\mathbb R}_q^{0|2}$, with an $\hat{R}$-matrix, can be expressed as follows
\begin{eqnarray} \label{49}
\sum_{k,l=1}^2 \hat{R}_{ij}^{kl} \, \theta_k \theta_l = q \theta_i \theta_j.
\end{eqnarray}

Using the $\hat{R}$-matrix, we can rewrite the relations (\ref{37}) and a relation in (\ref{39}) as follows
\begin{equation} \label{50}
\partial_{\theta_i} \theta_j = \delta_{ij} \, + q \sum_{k,l=1}^2 \hat{R}_{il}^{jk} \, \theta_l \partial_{\theta_k}, \quad
\partial_{\theta_i} \partial_{\theta_j} = q^2 \sum_{k,l=1}^2 \hat{R}_{ji}^{lk} \, \partial_{\theta_k} \partial_{\theta_l}.
\end{equation}

The proof of the following proposition is straightforward, so we omit it.
\begin{theorem}
$(a)$ $\hat{R}$ is invertible and satisfies the following equation (called the $Z_3$-graded braid relation):
\begin{equation*}
\hat{R}_{12}\hat{R}_{23}\hat{R}_{12} = \hat{R}_{23}\hat{R}_{12}\hat{R}_{23}.
\end{equation*}

\noindent$(b)$ The matrix $\hat{R}$ obeys the equation (called the $Z_3$-graded Hecke condition)
\begin{equation} \label{51}
(\hat{R}-qI)(\hat{R}+q^2I) = 0.
\end{equation}
\end{theorem}

The linear transformation $\hat{R}$ that acts on the square tensor space of the plane $\tilde{\mathbb R}_q^{0|2}$ with the properties in Theorem 7.1 is called a {\it Hecke symmetry}.

The Hecke condition implies that there is a vector space decomposition $V\otimes V\cong W_+\oplus W_-$ where $W_+$ and $W_-$ are eigenspaces for the eigenvalues $q$ and $-q^2$ of $\hat{R}$, respectively. In particular, $W_+ = $Im$(\hat{R}-qI)$ and $W_- = $Im$(\hat{R} + q^2I)$. An example of a  Hecke symmetry is the $Z_3$-graded permutation operator $\sigma: V\otimes V\longrightarrow V\otimes V$ defined by $\sigma(\theta_i\otimes \theta_j) = q^{\varrho(\theta_i)\varrho(\theta_j)} \theta_j\otimes \theta_i$ for all $i$ and $j$. Some extensions on these concepts will be considered in another study.

\section{A $Z_3$-graded universal enveloping algebra}\label{sec8}

In \cite{Jing}, a deformation of the universal enveloping algebra $U(gl(2))$ is constructed from the quantum
Weyl algebra.


In this section, we will construct a $Z_3$-graded universal enveloping algebra $U_q(\widetilde{gl}(2))$ using the quadratic elements of the Grassmann-Weyl algebra $\tilde{{\cal A}}_q(2)$. So, let us define the generators $E_{ij}$ for $i,j=1,2$ of $U_q(\widetilde{gl}(2))$ by the quadratic elements as follows
\begin{equation} \label{52}
E_{11} = \theta\partial_\theta, \quad E_{12} = \theta\partial_\varphi, \quad E_{21} = \varphi\partial_\theta, \quad E_{22} = \varphi\partial_\varphi.
\end{equation}
Then, we have
\begin{theorem} 
The defining relations for $U_q(\widetilde{gl}(2))$ as follows
\begin{eqnarray} \label{53}
&& E_{11}E_{12} - q^2 E_{12}E_{11} = E_{12}, \nonumber\\
&& E_{11}E_{21} - q E_{21}E_{11} = -qE_{21},  \nonumber\\
&& E_{11}E_{22} - E_{22}E_{11} = 0, \nonumber\\
&& E_{21}E_{12} - q^2 E_{12}E_{21} = E_{22} - E_{11} + (1-q^2) E_{11}^2, \nonumber \\
&& E_{12}E_{22} - E_{22}E_{12} = E_{12} + (q^2-1) E_{12}E_{11}, \nonumber \\
&& E_{21}E_{22} - E_{22}E_{21} = q E_{21} + (1-q) E_{21}E_{11}.
\end{eqnarray}
\end{theorem}

\noindent{\it Proof.} 
We only show the first one. Using the relations (\ref{37}) and (\ref{39}), one gets
\begin{eqnarray*}
E_{12}E_{11} &=& \theta\partial_\varphi \theta\partial_\theta = \theta^2 \partial_\theta \partial_\varphi\\
E_{11}E_{12}
&=& \theta\partial_\theta \theta\partial_\varphi = \theta(1+q^2 \, \theta \, \partial_\theta) \partial_\varphi = \theta\partial_\varphi + q^2 \theta^2
    \partial_\theta \partial_\varphi \\
&=& E_{12} + q^2 E_{12}E_{11}
\end{eqnarray*}
which is the first relation of (\ref{53}). \hfill $\square$

We see, from Theorem 5.1 and 5.2, that the element $D = E_{11}+E_{22} = \theta\partial_\theta + \varphi\partial_\varphi$ belongs to the center.

\end{document}